\documentclass{amsart}
\usepackage{amsmath,amscd,amsthm,amsfonts,amssymb}

\theoremstyle{plain}

\theoremstyle{definition}

\numberwithin{equation}{section}

\def \theo-intro#1#2 {\vskip .25cm\noindent{\bf Theorem #1\ }{\it #2}}

\def\fo#1{f \sb 0(#1)}

\def\g{\gamma}

\def\x{\xi}

\def\lpar{\left (}
\def\rpar{\right )}
\def\lmpar{\left \{ }
\def\rmpar{\right \}}
\def\lkpar{\left [}
\def\rkpar{\right ]}

\def\labs{\left |}
\def\rabs{\right |}
\def\babs#1{\labs {#1} \rabs}
\def\besfcn#1#2{J \sb {#1}({#2})}

\def\besker#1#2{\lpar 1 + {#1} \sp 2 \rpar \sp {#2}}

\def\dr{\, dr}

\def\drh{\, d\rho}

\def\dxi{\, d\x}
\def\dy{\, dy}

\def\dx{\, dx}

\def\R{{\mathbf R}}

\def\Rp{{\mathbf R \sb +}}

\def\ltworp{L \sp 2 \lpar \Rp \rpar}

\def\intbollone{\int \sb B}
\def\intr#1{\int \sb {\R \sp {#1}}}
\def\intrp{\int \sb 0 \sp \infty}

\def\coity{\mathcal C \sb 0 \sp \infty}

\def\sbsp#1#2#3{#2H \sp {#3}\lpar\R \sp {#1}\rpar}
\def\lsp#1#2{{L \sp {#2}\lpar\R \sp {#1}\rpar}}

\def\sop#1#2#3#4{#1S \sb {#2}#3#4}

\def\rop#1#2#3#4#5{[#1R \sb {#2} \sp {#3} {#4}\,]({#5})}
\def\nm#1#2#3{\left\|#1\right\| \sp {#2} \sb {#3}}

\def\ps{\psi}

\def\ch{\chi}

\def\r{\rho}

\def\rrh{r \rho}
\def\xt{[x](t)}

\def\xx{x \xi}
\def\z{\zeta}

\def\m{\mu}

\def\at{|t|}
\def\ax{|x|}
\def\axi{|\x|}

\def\l{\lambda}
\def\supp{\operatorname{supp}}

\begin{document}
\baselineskip 22pt \larger

\allowdisplaybreaks

\title
{Global Range Estimates for Maximal \\
 Oscillatory Integrals with Radial Testfunctions} \author
{Bj\"orn G.\ Walther \\ \\ \\
 Preprints in Mathematical Sciences 2011:1 ISSN 1403-9338 \\ \\
 LUNFMA-5068-2011} \date
{\today}  \address
{Centre for Mathematical Sciences, Lund University, SE--221 00 LUND, Sweden.} \email
{bjorn.walther@math.lu.se}

\begin{abstract}
We consider the maximal function of oscillatory integrals $S \sp a f$
where $\sop ( {} {\sp a f} ) \widehat \ (\x) = e \sp {it\axi \sp a} \widehat f(\xi)$ and $a \in \Rp \setminus 1$. For a fixed $n \geq 2$ we prove the global estimate
\begin{equation*}
\nm {S \sp a f} {} {L \sp 2 (\R \sp n, L \sp \infty (-1,1))} \, \leq \, C \nm f {} {\sbsp n {} s}, \quad s > a/4 \end{equation*}
with $C$ independent of the radial function $f$. We also prove that this result is almost sharp with respect to the Sobolev regularity $s$. \end{abstract} \maketitle \section
{Introduction} \subsection{}
In this paper we consider global range estimates for maximal functions of oscillatory integrals. More specifically we consider $\lsp n 2$-estimates of
\begin{equation}
\nm {\sop \lpar {} {\sp a f} \rpar [x]} {} {L \sp \infty (B)} \, = \, \sup \sb {\at < 1} \babs {\frac 1 {(2\pi) \sp n} \intr n e \sp {i(\xx + t\axi \sp a)} \widehat f(\x) \dxi}, \quad a > 0 \end{equation}
for $n \geq 2$. We will refer to $x$ as the {\it range} variable. $\widehat f$ is the Fourier transform of $f$,
\begin{equation}
\widehat f(\x) \, = \, \intr n e \sp {i\xx} f(x) \dx. \end{equation}
The test function $f$ will be radial. We obtain a linear estimate where the global range norm is controlled by a $\sbsp n {} s$-norm (Sobolev norm), and our result is almost sharp within the class of radial functions. \subsection {}
To obtain pointwise convergence results for oscillatory integrals of the type considered it is enough to consider local range estimates for maximal functions. However, global range estimates are of independent interest since they reveal global regularity properties of our oscillatory integrals. Global range estimates are also of interest when we for $a = 2$ consider the equivalence between local and global estimates due to Rogers \cite [theorem 3 p.\ 2108] {Rogers_2008}. We will suggest a generalisation of this result to other values of $a$. \subsection
{Earlier Results} The problem treated in this paper was for $a = 2$ introduced by Carleson \cite {Carleson_1980} and has been treated by many authors during the last couple of decades. See e.g.\ Ben-Artzi and Devinatz \cite {Ben-Artzi_Devinatz_1991}, Bourgain \cite {Bourgain_1992}, Cho, Lee and Shim \cite {Cho_Lee_Shim_2006}, Dahlberg and Kenig \cite {Dahlberg_Kenig_1982}, Gigante and Soria \cite {Gigante_Soria_2008}, Kenig, Ponce and Vega \cite {Kenig_Ponce_Vega_1991}, Kolasa \cite {Kolasa_1997}, Lee \cite {Lee_2006}, Moyua, Vargas and Vega \cite {Moyua_Vargas_Vega_1996}, Prestini \cite {Prestini_1990}, Rogers \cite {Rogers_2008}, Rogers and Villarroya \cite {Rogers_Villarroya_2008}, Sj\"olin \cite {Sjoelin_1987, Sjoelin_1997, Sjoelin_2005}, Tao and Vargas \cite {Tao_Vargas_2000}, \cite {Walther_1999b, Walther_2001}, S.\ Wang \cite {Wang_Sichun_2006} and the papers cited there. In some of these papers (e.g.\ \cite {Kenig_Ponce_Vega_1991}, \cite {Rogers_Villarroya_2008}, \cite {Sjoelin_1997} and \cite {Walther_2001}) $L \sp q$-estimates are considered for some $q \neq 2$. We will however restrict ourselves to the case $q = 2$. \par
Estimates which are sharp or almost sharp with respect to the number of derivatives $s$ have been obtained in the cases $n = 1$ and $a = 1$ only. See \S\S \ref {subsubsec_1.3.3} and \ref {subsubsec_1.3.4} for a discussion on these estimates. Some of these results will be used as tools in our proofs. \subsubsection
{Some best known local range results} The best known local range result for $a = n = 2$ is due to Lee \cite {Lee_2006} where it is proved that $s > 3/8$ is sufficient for a norm inequality. The best known local range result for fixed $n \geq 2$ and for fixed $a > 1$ is that $s > 1/2$ is sufficient for a norm inequality. This result is due to Sj\"olin \cite {Sjoelin_1987} and Vega \cite {Vega_1988}. See also S.L.\ Wang \cite {Wang_Si_Lei_1991}, \cite {Walther_1999b} and \cite [Example 4.1, p.\ 331] {Walther_2002b}. Soljanik has proved a local range estimate with $s > 1/2$ as sufficient condition where $m(t,\r) = \exp(it\r \sp a)$ may be replaced by a function $m$ which is assumed to be bounded only. See \cite [theorem 14.3 p.\ 219] {Walther_1999a}. \subsubsection
{Some best known global range results} The best known global range result for $a = n = 2$ is that $s > 3/4$ is sufficient for a norm inequality. This follows from the local range result of Lee \cite {Lee_2006} and from the equivalence result of Rogers cited as theorem \ref {thm_3.02} on page \pageref {thm_3.02}. The best known global range result for fixed $n \geq 2$ and for fixed $a > 0$ is that $s > a/2$ is sufficient for a norm inequality. See Carbery \cite {Carbery_1985} and the references cited in theorem \ref {thm_4.02} on page \pageref {thm_4.02}. \subsubsection
{Sharp and almost sharp results for $n = 1$.} \label {subsubsec_1.3.3} Consider first the case $a < 1$. If $s > a/4$ then $\nm {\sop ( {} {\sp a f} ) [x]} {} {L \sp \infty (B)}$ can be estimated globally. See theorem \ref {thm_4.04} on page \pageref {thm_4.04}. (Recall that we consider $L \sp 2$-estimates only.) On the other hand, if $s < a/4$ then $\nm {\sop ( {} {\sp a f} ) [x]} {} {L \sp \infty (B)}$ cannot be estimated even locally. See \cite [theorem 1.2 (b) p.\ 486] {Walther_1995}. The interval $]a/4, \infty[$ is thus the largest open interval of admissable regularities $s$ in the local as well as in the global case. When $s = a/4$ the existence of an estimate is an open problem in the local as well as in the global case. \par
Consider now the case $a > 1$. If $s = 1/4$ then $\nm {\sop ( {} {\sp a f} ) [x]} {} {L \sp \infty (B)}$ can be estimated locally. See Sj\"olin \cite [theorem 3 p.\ 700] {Sjoelin_1987}. If $s > a/4$ then $\nm {\sop ( {} {\sp a f} ) [x]} {} {L \sp \infty (B)}$ can be estimated globally. See theorem \ref {thm_4.04} on page \pageref {thm_4.04}. On the other hand, if $s < 1/4$ then $\nm {\sop ( {} {\sp a f} ) [x]} {} {L \sp \infty (B)}$ cannot be estimated locally (\cite [theorem 4 p.\ 700] {Sjoelin_1987}), and if $s < a/4$ then $\nm {\sop ( {} {\sp a f} ) [x]} {} {L \sp \infty (B)}$ cannot be estimated globally (\cite [p.\ 106] {Sjoelin_1994}). Thus the interval $]1/4, \infty[$ is the largest open interval of admissable regularities in the local case, and the interval $]a/4, \infty[$ is the largest open interval of admissable regularities in the global case. When $s = a/4$ the existence of an estimate is an open problem in the global case. \subsubsection
{Sharp results for $a = 1$.} \label {subsubsec_1.3.4} Fix $n \geq 1$. If $s > 1/2$ then $\nm {\sop ( {} {\sp a f} ) [x]} {} {L \sp \infty (B)}$ can be estimated globally. See theorem \ref {thm_4.02} on page \pageref {thm_4.02}. This result is well known. On the other hand, if $s = 1/2$ then $\nm {\sop ( {} {\sp a f} ) [x]} {} {L \sp \infty (B)}$ cannot be estimated even locally. See \cite [theorem 14.2 p.\ 216] {Walther_1999a}. Thus the interval $]1/2, \infty[$ is the largest open interval of admissable regularities in the local as well as in the global case. \subsection
{The Plan of this Paper} In \S \ref {sec_2} we introduce notation and state our results. The problem we study is in part motivated by the equivalence between local and global range estimates in a special case. A brief discussion on this equivalence is found in \S \ref {sec_3}. In \S \ref {sec_4} we have collected results needed in our proofs and in \S \ref {sec_5} we prove our result. \section
{Notation and Results} \label {sec_2} \subsection {}
In this section we introduce some notation used in this paper and formulate our result which is almost sharp within the class of radial functions. \par
Unless otherwise explicitly stated all functions $f$ and $g$ are supposed to belong to $\coity (\R \sp n \setminus 0)$. \subsection
{Oscillatory Integrals, The Fourier Transform and Inhomogeneous Sobolev Spaces} For the range variable $x \in \R \sp n$ and $t \in \R$ we define
\begin{align}
\sop ( {} {\sp a f} ) \xt \, &= \, \frac 1 {(2\pi) \sp n} \intr n e \sp {i(\xx + t\axi \sp a)} \widehat f(\x) \dxi \intertext
{Here $\widehat f$ is the Fourier transform of $f$,}
\widehat f(\x) \, &= \, \intr n e \sp {-i\xx} f(x) \dx. \end{align}
We also introduce inhomogeneous fractional $\lsp n 2$-based Sobolev spaces
\begin{equation}
\sbsp n {} s \, = \, \lmpar f \in \mathcal S \sp \prime (\R \sp n) : \intr n \besker \axi s \babs {\widehat f(\xi)} \sp 2 \dxi < \infty \rmpar.
\end{equation} \subsection
{Auxiliary Notation} $B \sp n$ denotes the open unit ball in $\R \sp n$. We write $B \sp 1 = B$. Throughout this paper we will use auxiliary functions $\ch$ and $\ps$ such that $\ch \in \coity (\R)$ is even,
\begin{equation}
\chi(\R \setminus 2B) \, = \, \{ 0 \}, \quad \chi(\R) \, \subseteq \, [0,1], \quad \chi(B) \, = \, \{ 1 \} \end{equation}
and $\ps = 1 - \ch$. From $\ch$ we obtain a family of functions as follows: for a fixed $m > 1$ set $\ch \sb m (\x) = \ch (\x/m)$. \par
Numbers denoted by $C$ (sometimes with subscripts) may be different at each occurrence. \subsection
{Theorem A} \label {thm_2.04} {\it Let $a \in \R \sb + \setminus 1$, $n \geq 2$ and $s > a/4$. Then there is a number $C$ independent of $f$ in the class of radial functions such that
\begin{equation}
\nm {S \sp a f} {} {L \sp 2 (\R \sp n, L \sp \infty (B))} \, \leq \, C \nm f {} {\sbsp n {} s}. \end{equation}} \subsection
{Remark} \label {subsec_2.5} Fix $a > 1$. If $s > a/4$ then
\begin{equation}
q \, > \, \frac {4(a - 1)n} {4s + a(2n - 1) - 2n} \, = \, \frac {4n(a - 1)} {4s - a + 2n(a - 1)} \end{equation}
for $q = 2$. Hence theorem A in the case $a > 1$ follows from the sufficiency part of Sj\"olin \cite [theorem 4 p.\ 37] {Sjoelin_1997}. We have chosen to include the case $a > 1$ here since the proof uses a reduction to the case $n = 1$, and hence we can use theorem \ref {thm_4.04} on page \pageref {thm_4.04}. This reduction follows the same pattern regardless of the choice of $a$. \subsection
{Theorem B} \label {thm_2.06} {\it Let $a \in \R \sb + \setminus 1$ and $n \geq 2$. Assume that there is a number $C$ independent of $f$ in the class of radial functions such that
\begin{equation} \label {eq_2.06}
\nm {S \sp a f} {} {L \sp 2 (\R \sp n, L \sp \infty (B))} \, \leq \, C \nm f {} {\sbsp n {} s}. \end{equation}
Then $s \geq a/4$.} \subsection
{Remark} \label {subsec_2.07} Consider again the case $a > 1$. If
\begin{align}
q \, &< \, \frac {4(a - 1)n} {4s + a(2n - 1) - 2n} \intertext
{for $q = 2$, i.e. if}
\label {eq_2.09}
2 \, &< \, \frac {4n(a - 1)} {4s - a + 2n(a - 1)} \end{align}
then $-2n(a - 1) < 4s - a < 0$. Conversely, if $-2n(a - 1) < 4s - a < 0$ then \eqref {eq_2.09} holds. Hence theorem B in the case $a > 1$ follows from the necessity part of Sj\"olin \cite [theorem 4 p.\ 37] {Sjoelin_1997}. \section
{A Brief Discussion on Equivalence between \\
 Local and Global Range Estimates} \label {sec_3} \subsection {}
Let us say that the open interval $I$ is {\it locally admissable} if $s \in I$ implies that there is a number $C$ independent of $f$ such that
\begin{equation}
\nm {S \sp a f} {} {L \sp 2 (B \sp n, L \sp \infty (B))} \, \leq \, C \nm f {} {\sbsp n {} s}. \end{equation}
If instead this estimate holds with $B \sp n$ replaced by $\R \sp n$ then we say that $I$ is {\it globally admissable}. With this terminology we have e.g.\ that $I = ]a/4, \infty[$ is locally and globally admissable when $a < 1 = n$. Moreover, $I$ is maximal with this property. See \S \ref {subsubsec_1.3.3} on page \pageref {subsubsec_1.3.3}. \subsection
{Theorem \rm (Rogers \cite [theorem 3 p.\ 2108] {Rogers_2008})} \label {thm_3.02} {\it Let $a = 2$ and let $n \geq 1$ be fixed. Then the interval $]\sigma, \infty[$ is locally admissable if and only if the interval $]2\sigma, \infty[$ is globally admissable.} \subsection {}
Note that theorem \ref {thm_3.02} is consistent with results explained in \S \ref {subsubsec_1.3.3} on page \pageref {subsubsec_1.3.3} where {\it maximal} admissable intervals are given. One may conjecture that for a fixed $a > 1$ and for a fixed $n \geq 1$ the interval $]\sigma, \infty[$ is locally admissable if and only if the interval $]a\sigma, \infty[$ is globally admissable. For $n = 1$ that conjecture holds true (see \S \ref {subsubsec_1.3.3} on page \pageref {subsubsec_1.3.3}) and is consistent with the following conjectures which hold true for the subclass of radial test functions (cf.\ theorem A on page \pageref {thm_2.04} and theorem \ref {thm_4.07} on page \pageref {thm_4.07}): \subsection
{Conjecture 1} {\it Let $a < 1$. Then the interval $]a/4, \infty[$ is globally admissable.} \subsection
{Conjecture 2} {\it Let $a > 1$. Then the interval $]1/4, \infty[$ is locally admissable, and the interval $]a/4, \infty[$ is globally admissable.} \subsection {}
We end this brief discussion on equivalence between local and global range estimates by noting that the maximal intervals which are locally and globally admissable coincide if $a < 1 = n$ or if $a = 1$. \section
{Preparation} \label {sec_4} \subsection {}
In this section we collect results needed to give proofs of the theorems stated in section 
\ref {sec_2}. \subsection
{Theorem \rm (Cf.\ Cowling \cite {Cowling_1983}, Cowling, Mauceri \cite {Cowling_Mauceri_1985}, Rubio de Francia \cite {Rubio_de_Francia_1986}, Sogge, Stein \cite {Sogge_Stein_1985}, Stein \cite [\S XI.4.1, p.\ 511] {Stein_1993} and \cite [theorem 14.1 p.\ 215] {Walther_1999a}.)} \label {thm_4.02} {\it Assume that the functions $w \sb 1$ and $w \sb 2$ belong to $\lsp {} 2$ and that the function $m$ satisfies the following assumption: there is a number $C$ independent of $(t,\x)$ such that
\begin{equation}
|m(t,\x)| \, \leq \, Cw \sb 1 (t), \quad \babs{\lkpar \partial \sb 1 m \rkpar (t,\x)} \, \leq \, C \lpar w \sb 1 (t) + w \sb 2 (t) \axi \sp a \rpar, \quad a > 0. \end{equation}
If $s > a/2$, then there is a number $C$ independent of $f$ such that
\begin{equation}
\lpar \intr n \sup \sb {t \in B} \babs {\intr n e \sp {i\xx} m(t,\xi) \widehat f(\x) \dxi} \sp 2 \dx \rpar \sp {1/2} \, \leq \, C \nm f {} {\sbsp n {} s}. \end{equation}} \subsection
{Corollary} \label {cor_4.03} {\it Assume that $m$ fullfills the same assumptions as in theorem \ref {thm_4.02}. Then there is a number $C$ independent of $f$ such that
\begin{equation}
\lpar \intr n \sup \sb {t \in B} \babs {\intr n e \sp {i\xx} m(t,\xi) \widehat f(\x) \dxi} \sp 2 dx \rpar \sp {1/2} \leq C \nm f {} {\lsp n 2}, \quad \supp \widehat f \subseteq 2B \sp n. \end{equation}} \subsection*
{\it Proof:} If $\supp \widehat f \subseteq 2B \sp n$, then $f \in \sbsp n {} s$ for every $s$. We now apply theorem \ref {thm_4.02} in conjunction with the estimate
\begin{equation}
\nm f {} {\sbsp n {} s} \, \leq \, C \nm f {} {\lsp n 2}, \quad \supp \widehat f \subseteq 2B \sp n \end{equation}
where the number $C$ is independent of $f$. \subsection
{Theorem \rm (Cf.\ \cite [theorem 2.5 p.\ 159] {Walther_2002a} for $a < 1$ and Sj\"olin \cite [p.\ 106] {Sjoelin_1994} for $a > 1$.)} \label {thm_4.04} {\it Let $a \in \R \sb + \setminus 1$ and $s > a/4$. Then there is a number $C$ independent of $f$ such that
\begin{equation}
\nm {S \sp a f} {} {L \sp 2 (\R, L \sp \infty (B))} \, \leq \, C \nm f {} {\sbsp {}{} s}. \end{equation}} \subsection
{\it Remarks on the Proof:} In the case $a < 1$ the proof is found in \cite [\S 4 pp.\ 161--164] {Walther_2002a} and in the case $a > 1$ the proof is found in Sj\"olin \cite [pp.\ 107--112] {Sjoelin_1994}. An important tool in both cases is the smooth decomposition of Littlewood and Paley. More precisely, if $N$ is a dyadic integer and $\eta \in \coity (\R)$ is an even function such that
\begin{equation}
\eta \lpar \R \sb + \setminus [1/2,2] \rpar \, = \, \{ 0 \}, \quad \eta(\R) \, \subseteq \, [0,1] \end{equation}
and
\begin{equation}
\sum \sb {N > 1} \eta(N\x) + \sum \sb N \eta(\x/N) \, = \, 1, \quad \x \neq 0 \end{equation}
(cf.\ Bergh, L\"ofstr\"om \cite [Lemma 6.1.7, pp.\ 135--136] {Bergh_Loefstroem}) then there are positive numbers $C \sb 1$ and $C \sb 2$ independent of $f$ such that
\begin{equation}
C \sb 1 \nm f 2 {\sbsp {}{} s} \, \leq \, \intr {} \lpar \ch(\x) + \sum \sb N N \sp {2s} \eta(\x/N) \rpar \babs {\widehat f(\x)} \sp 2 \dx \, \leq \, C \sb 2 \nm f 2 {\sbsp {}{} s}. \end{equation}
We denote the expression within the brackets by $\g \sb {2s} (\x)$. Further elaboration (Parseval's formula, approximation of operators using cutoff functions, passing to the adjoint, Fatou's lemma, and Fubini's theorem) leads to deriving a $m\m$-uniform $\lsp {} 1$-estimate for the kernel
\begin{equation}
K \sb {m\m}(x) \, = \, \ch \sb m (x) \sup \sb {t \in 2B} \babs {\intr {} e \sp {i\xx} e \sp {it\axi \sp a} \g \sb {-2s} (\x) \ch \sb \m (\x) \sp 2 \dxi}. \end{equation}
Here the low and high frequency contributions are analyzed separately. Considering $\g \sb {-2s}$, the low frequency contribution corresponds to the term $\ch$ and the high frequency contribution to the infinite sum over the dyadic integers. \subsubsection {}
It is useful to note that we may avoid analyzing the low frequency contribution in the way done in \cite [Lemma 3.2 p.\ 160] {Walther_2002a} and Sj\"olin \cite [p.\ 109--110] {Sjoelin_1994}. More precisely, after having linearised the maximal operator so as to obtain the operator
\begin{equation}
\rop {} t {} f x \, = \, \intr {} e \sp {i\xx} e \sp {it(x)\axi \sp a} \g \sb {-2s} (\x) \sp {1/2} f(\x) \dxi \end{equation}
we may define
\begin{equation}
R \sb {t, \z} f \, = \, R \sb t (\z f), \quad \z \in \{ \ch, \ps \}. \end{equation}
(As usual, $t$ is any measurable function such that $t(\R) \subseteq B$ and we want to find estimates independent of $t$.) Then we use corollary \ref {cor_4.03} on page \pageref {cor_4.03} to estimate $R \sb {t, \ch} f$. The high frequency contribution in $K \sb {m\m}$ is used to estimate $R \sb {t, \ps}f$. \par
The use of $R \sb {t, \z}$ for a fixed $\z \in \{ \ch, \ps \}$ of course corresponds to the decomposition
\begin{equation}
\widehat f \, = \, \ch \widehat f + \ps \widehat f. \end{equation} \subsection
{Theorem \rm (\cite [theorem C p.\ 190] {Walther_2001})} \label {thm_4.06} {\it Let $a < 1$, $n \geq 2$ and $s > a/4$. Then there is a number $C$ independent of $f$ in the class of radial functions such that
\begin{equation}
\nm {S \sp a f} {} {L \sp 2 (B \sp n, L \sp \infty (B))} \, \leq \, C \nm f {} {\sbsp n {} s}. \end{equation}} \subsection
{Theorem \rm (Sj\"olin \cite [theorem 1 p.\ 135] {Sjoelin_1995})} \label {thm_4.07} {\it Let $a > 1$ and $n \geq 2$. Then there is a number $C$ independent of $f$ in the class of radial functions such that
\begin{equation}
\nm {S \sp a f} {} {L \sp 2 (B \sp n, L \sp \infty (B))} \, \leq \, C \nm f {} {\sbsp n {} {1/4}}. \end{equation}} \subsection
{Theorem \rm (Cf.\ e.g.\ Stein, Weiss \cite [theorem 3.10 p.\ 158] {Stein_Weiss}.)} \label {thm_4.08} {\it Let $f$ be radial. Then
\begin{equation}
\widehat f(\x) \, = \, (2\pi) \sp {n/2} \axi \sp {-n/2 + 1} \intrp \fo r \besfcn {n/2 - 1} {r\axi} r \sp {n/2} \dr \end{equation}
where $f(x) = \fo {|x|}$ and $J \sb \l$ is the Bessel function of the first kind of order $\l$.} \subsection
{Theorem \rm (Cf.\ e.g.\ Stein, Weiss \cite [Lemma 3.11 p.\ 158] {Stein_Weiss}.)} \label {thm_4.09} {\it If $\l > -1/2$, then there is a number $C \sb \l$ independent of $\r > 1$ such that
\begin{equation}
\babs {\besfcn \l \r - \lpar \frac 2 \pi \rpar \sp {1/2} \r \sp {-1/2} \cos \lpar \r - \frac {\l\pi} 2 - \frac \pi 4 \rpar} \, \leq \, C \sb \l \r \sp {-3/2}. \end{equation}} \subsection
{Theorem \rm (\cite [theorem 2.6 p.\ 3644] {Walther_2002c})} \label {thm_4.10} {\it Define $f \sb y$ by
\begin{equation}
\widehat {f \sb y}(\x) \, = \, e \sp {iy\axi} \widehat f(\x), \quad y \in B.
\end{equation}
Assume that $n \geq 2$, $a < 1$ and that there is a number $C$ independent of the radial function $f$ such that
\begin{equation}
\intbollone \nm {S \sp a f \sb y} 2 {L \sp 2 (B \sp n, L \sp \infty (B))} \dy \, \leq \, C \, \nm f 2 {\sbsp n {} s}. \end{equation}
Then $s \geq a/4$.} \section
{Proofs} \label {sec_5} \subsection {\it
 Proof of theorem A in \S {\rm \ref {thm_2.04}} on page \pageref {thm_2.04}} Assume that $s > a/4$. Define
\begin{equation}
\lpar \widetilde {S \sp a} f \rpar \xt \, = \, \intr n e \sp {i(\xx + t\ax \sp a)} \besker \x {-s/2} f(\x) \dxi. \end{equation}
To prove the theorem it is according to Parseval's formula enough to prove that there is a number $C$ independent of the radial function $f$ such that
\begin{equation} \label {eq_5.02}
\nm {\widetilde {S \sp a} f} {} {L \sp 2 (\R \sp n, L \sp \infty (B))} \, \leq \, C \nm f {} {\lsp n 2}. \end{equation}
Define
\begin{equation}
\lpar \widetilde {S \sb \z \sp a} f \rpar \xt \, = \, \intr n e \sp {i(\xx + t\ax \sp a)} \besker \x {-s/2} \z (\axi) f(\x) \dxi, \quad \z \in \{ \ch, \ps \}. \end{equation}
It is then sufficient to prove the estimate \eqref {eq_5.02} with $\widetilde {S \sp a}$ replaced by $\widetilde {S \sb \z \sp a}$ for all $\z \in \{ \ch, \ps \}$.

The estimate for $\widetilde {S \sb \ch \sp a}$ follows from corollary \ref {cor_4.03} on page \pageref {cor_4.03}. Hence it remains to prove the estimate for $\widetilde {S \sb \ps \sp a}$.

Define
\begin{equation}
\lpar \widetilde {{} \sb \z S \sb \ps \sp a} f \rpar [x] \, = \, \z(\ax) \lpar \widetilde {S \sb \ps \sp a} f \rpar [x]. \end{equation}
It is then sufficient to prove the estimate \eqref {eq_5.02} with $\widetilde {S \sp a}$ replaced by $\widetilde {{} \sb \z S \sb \ps \sp a}$ for all $\z \in \{ \ch, \ps \}$.

The estimate for $\widetilde {{} \sb \ch S \sb \ps \sp a}$ follows from theorem \ref {thm_4.06} on page \pageref {thm_4.06} in the case $a < 1$ and from theorem \ref {thm_4.07} on page \pageref {thm_4.06} in the case $a > 1$. Hence it remains to prove the estimate for $\widetilde {{} \sb \ps S \sb \ps \sp a}$. \par
There is a function $f \sb 0 \in \coity (\Rp)$ such that
\begin{equation}
f(\x) \, = \, \axi \sp {-n/2 + 1/2} \fo {\axi}. \end{equation}
According to theorem \ref {thm_4.08} on page \pageref {thm_4.08}
\begin{multline}
\lpar \widetilde {{} \sb \ps S \sb \ps \sp a} f \rpar \xt \, = \\
= (2\pi) \sp {n/2} r \sp {-n/2 + 1} \ps(r) \intrp \r \sp {1/2} \besfcn {n/2 - 1} \rrh e \sp {it\r \sp a} \besker \r {-s/2} \ps(\r) \fo \r \drh, \end{multline}
where $r = \ax$. Let $t$ be any measurable function such that $t(\Rp) \subseteq B$ and define
\begin{equation}
\rop {} t {} f r \, = \, \intrp \ps(r) (\rrh) \sp {1/2} \besfcn {n/2 - 1} \rrh e \sp {it(r)\r \sp a} \r \sp {-s} \ps(\r) f(\r) \drh \end{equation}
for $f \in \coity (\Rp)$. To prove the estimate for $\widetilde {{} \sb \ps S \sb \ps \sp a}$ it is sufficient to prove that there is a number $C$ independent of $f$ {\it and} $t$ such that
\begin{equation}
\nm {R \sb t f} {} \ltworp \, \leq \, C \nm f {} \ltworp. \end{equation}
Define
\begin{equation}
\rop {} {t,1} {} f r \, = \, \lpar \frac 2 \pi \rpar \sp {1/2} \intrp \ps(r) \cos \lpar \rrh - \frac {\l \pi} 2 - \frac \pi 4 \rpar e \sp {it(r)\r \sp a} \r \sp {-s} \ps(\r) f(\r) \drh \end{equation}
and $R \sb {t,2} = R \sb t - R \sb {t,1}$ where $f \in \coity (\Rp)$. \subsubsection {} \label {subsubsec_5.1.1}
Due to the cosine factor there are numbers $C \sb 1$ and $C \sb 2$ independent of $r$, $f$ and $t$ such that
\begin{multline*}
\rop {} {t,1} {} f r \, = \, C \sb 1 \intrp \ps(r) e \sp {i(\rrh + t(r)\r \sp a)} \r \sp {-s} \ps(\r) f(\r) \drh \, + \\
+ \, C \sb 2 \intrp \ps(r) e \sp {i(-\rrh + t(r)\r \sp a)} \r \sp {-s} \ps(\r) f(\r) \drh. \end{multline*}
Here we apply theorem \ref {thm_4.04} on page \pageref {thm_4.04} to each term in the right hand side. We get that there is a number $C$ independent of $f$ and $t$ such that
\begin{equation}
\nm {R \sb {t,1} f} {} \ltworp \, \leq \, C \nm f {} \ltworp. \end{equation} \subsubsection {}
It only remains to prove that there is a number $C$ independent of $f$ and $t$ such that
\begin{equation}
\nm {R \sb {t,2} f} {} \ltworp \, \leq \, C \nm f {} \ltworp. \end{equation}
According to theorem \ref {thm_4.09} on page \pageref {thm_4.09} there is a number $C$ independent of $f$ such that
\begin{align}
\babs {\rop {} {t,2} {} f r} \, &\leq \, \frac {C \ps(r)} r \intrp \r \sp {-1 - s} \babs{f(\r)} \ps(\r) \drh. \intertext
{We use Cauchy-Schwarz inequality to get}
\babs {\rop {} {t,2} {} f r} \, &\leq \, \frac {C \ps(r)} r \lpar \intrp \r \sp {-2 - 2s} \ps(\r) \drh \rpar \sp {1/2} \nm f {} {\ltworp}. \intertext
{Squaring and integrating with respect to $r$ gives that there is a number $C$ independent of $f$ such that}
\nm {R \sb {t,2} f} 2 {\ltworp} \, &\leq \, C \lpar \intrp \frac {\ps(r) \sp 2 \dr} {r \sp 2} \rpar \nm f 2 {\ltworp}. \end{align} \subsection {\it 
Proof of theorem B in \S {\rm \ref {thm_2.06}} on page \pageref {thm_2.06} in the case $a < 1$.} Assume that there is a number $C$ independent of $f$ in the class of radial functions such that
\begin{align}
\nm {S \sp a f} {} {L \sp 2 (\R \sp n, L \sp \infty (B))} \, &\leq \, C \nm f {} {\sbsp n {} s}. \intertext
{In particular we have}
\label {eq_5.16}
\nm {S \sp a f \sb y} {} {L \sp 2 (\R \sp n, L \sp \infty (B))} \, &\leq \, C \nm {f \sb y} {} {\sbsp n {} s} \, = \, C \nm f {} {\sbsp n {} s} \end{align}
where
\begin{equation}
\widehat {f \sb y}(\x) \, = \, e \sp {iy\axi} \widehat f(\x), \quad y \in B. \end{equation}
Squaring \eqref {eq_5.16} and integrating with respect to $y$ gives that there is a number $C$ independent of $f$ such that
\begin{equation}
\intbollone \nm {S \sp a f \sb y} 2 {L \sp 2 (\R \sp n, L \sp \infty (B))} \dy \, \leq \, C \nm f 2 {\sbsp n {} s}, \end{equation}
and the inequality remains valid when we replace $\R \sp n$ in the left hand side by $B \sp n$. The conclusion sought for now follows from theorem \ref {thm_4.10} on page \pageref {thm_4.10}. \subsection {\it 
Proof of theorem B in \S {\rm \ref {thm_2.06}} on page \pageref {thm_2.06} in the case $a > 1$} See remark \ref {subsec_2.07} on page \pageref {subsec_2.07}.

\begin{thebibliography}{99}

\bibitem
{Ben-Artzi_Devinatz_1991}
 M. Ben-Artzi, A. Devinatz, \textit
{Local Smoothing and Convergence Properties of Schrödinger Type Equations,}
 J. Funct. Anal. \textbf
{101} (1991),
 231--254, MR \textbf
{92k}:35064.

\bibitem
{Bergh_Loefstroem}
 J. Bergh, J. L\"ofstr\"om, \textit
{Interpolation theory. An Introduction,} Grundlehren der Mathematischen Wissenschaften, No. 223. Springer-Verlag, Berlin-New York, 1976, MR \textbf
{58} \#2349.

\bibitem
{Bourgain_1992}
 J. Bourgain, \textit
{A Remark on Schr\"odinger Operators,}
 Israel J. Math. \textbf
{77} (1992),
 1--16, MR \textbf
{93k}:35071.

\bibitem
{Carbery_1985}
 A. Carbery, \textit
{Radial Fourier Multipliers and Associated Maximal Functions,}
 Recent progress in Fourier analysis (El Escorial, 1983), 49--56,
 North-Holland Math. Stud., \textbf
{111}, North-Holland, Amsterdam, 1985, MR \textbf
{87i}:42029.

\bibitem
{Carleson_1980}
 L. Carleson, \textit
{Some Analytic Problems Related to Statistical Mechanics,}
 Euclidean Harmonic Analysis (Proc. Sem., Univ. Maryland, College Park, Md., 1979), 5--45,
 Lecture Notes in Math. \textbf
{779}, Springer, Berlin, 1980, MR \textbf
{82j}:82005.

\bibitem
{Cho_Lee_Shim_2006}
 Y. Cho, S. Lee, Y. Shim, \textit
{A Maximal Inequality Associated to Schr\"odinger Type Equation,}
 Hokkaido Math. J. \textbf
{35} (2006),
 767-–778, MR \textbf
{2007m}:42016.

\bibitem
{Cowling_1983}
 M. Cowling, \textit
{Pointwise Behavior of Solutions to Schr\"odinger Equations,}
 Harmonic Analysis (Cortona, 1982), 83--90,
 Lecture Notes in Math. \textbf
{992}, Springer, Berlin, 1983, MR \textbf
{85c}:34029.

\bibitem
{Cowling_Mauceri_1985}
 M. Cowling, G. Mauceri, \textit
{Inequalities for some Maximal Functions. I,}
 Trans. Amer. Math. Soc. \textbf
{287} (1985), 
 431--455, MR \textbf
{86a}:42023.

\bibitem
{Dahlberg_Kenig_1982}
 B. Dahlberg, C.E. Kenig, \textit
{A Note on the Almost Everywhere Behavior of Solutions to the Schrödinger Equation,}
 Harmonic Analysis (Minneapolis, Minn., 1981), 205--209,
 Lecture Notes in Math. \textbf
{908}, Springer, Berlin-New York, 1982, MR \textbf
{83f}:35023.

\bibitem
{Gigante_Soria_2008}
 G. Gigante, F. Soria, \textit
{On the Boundedness in $H \sp {1/4}$ of the Maximal Square Function Associated with the Schr\"odinger Equation,}
 J. Lond. Math. Soc. (2) \textbf
{77} (2008), 
 51--–68, MR \textbf
{2010a}:42039.
 
\bibitem
{Kenig_Ponce_Vega_1991}
 C.E. Kenig, G. Ponce, L. Vega, \textit
{Oscillatory Integrals and Regularity of Dispersive Equations,}
 Indiana Univ. Math. J. \textbf
{40} (1991), 
 31--69, MR \textbf
{98h}:35043.

\bibitem
{Kolasa_1997}
 L. Kolasa, \textit
{Oscillatory Integrals and Schr\"odinger Maximal Operators,}
 Pacific J. Math. \textbf
{177} (1997), 
 77--101, MR \textbf
{98h}:35043.

\bibitem
{Moyua_Vargas_Vega_1996}
 A. Moyua, A. Vargas, L. Vega, \textit
{Schr\"odinger Maximal Function and Restriction Properties of the Fourier Transform,}
 Internat. Math. Res. Notices \textbf
{1996},
 793--815, MR \textbf
{97k}:42042.

\bibitem
{Lee_2006}
 Lee, S., \textit
{On Pointwise Convergence of the Solutions to Schr\"odinger Equations in $\R \sp 2$,}
 Int. Math. Res. Not. \textbf
{2006}, Article ID 32597, 21 pp., MR \textbf
{2007j}:35180.

\bibitem
{Prestini_1990}
 E. Prestini, \textit
{Radial Functions and Regularity of Solutions to the Schrödinger Equation,}
 Monatsh. Math. \textbf
{109} (1990),
 135--143, MR \textbf
{91j}:35035.

\bibitem
{Rubio_de_Francia_1986}
 J.L. Rubio de Francia, \textit
{Maximal Functions and Fourier Transforms,}
 Duke Math. J. \textbf
{53} (1986),
 395--404, MR \textbf
{87j}:42046.

\bibitem
{Rogers_2008}
 K.M. Rogers, \textit
{A Local Smoothing Estimate for the Schr\"odinger Equation,}
 Adv. Math. \textbf
{219} (2008),
 2105--2122, MR \textbf
{2009j}:35055.

\bibitem
{Rogers_Villarroya_2008}
 K.M. Rogers, P. Villarroya \textit
{Sharp Estimates for Maximal Operators Associated to the Wave Equation,}
 Ark. Mat. \textbf
{46} (2008),
 143--151, MR \textbf
{2009b}:35242.

\bibitem
{Sjoelin_1987}
 P. Sj\"olin, \textit
{Regularity of Solutions to the Schr\"odinger Equation,}
 Duke Math. J. \textbf
{55} (1987),
 699--715, MR \textbf
{88j}:35026.

\bibitem
{Sjoelin_1994}
 P. Sj\"olin, \textit
{Global Maximal Estimates for Solutions to the Schr\"odinger Equation,}
 Studia Math. \textbf
{110} (1994),
 105--114, MR \textbf
{95e}:35052.

\bibitem
{Sjoelin_1995}
 P. Sj\"olin, \textit
{Radial Functions and Maximal Estimates for Solutions to the Schr\"odinger Equation,}
 J. Austral. Math. Soc. Ser. A \textbf
{59} (1995),
 134--142, MR \textbf
{96d}:42032.

\bibitem
{Sjoelin_1997}
 P. Sj\"olin, \textit
{$L \sp p$ Maximal Estimates for Solutions to the Schr\"odinger Equation,}
 Math. Scand. \textbf
{81} (1997),
 35--68, MR \textbf
{98j}:35038.

\bibitem
{Sjoelin_2005}
 P. Sj\"olin, \textit
{Spherical Harmonics and Maximal Estimates for the Schr\"odinger Equation,}
 Ann. Acad. Sci. Fenn. Math. \textbf
{30} (2005),
 393--406, MR \textbf
{2006g}:42035.

\bibitem
{Sogge_Stein_1985}
 C.D. Sogge, E.M. Stein, \textit
{Averages of Functions over Hypersurfaces of $\R \sp n$,}
 Invent. Math. \textbf
{82} (1985),
 543--556, MR \textbf
{87d}:42030.

\bibitem
{Stein_1993} E. M. Stein, \textit
{Harmonic Analysis: Real-variable methods, Orthogonality, and Oscillatory integrals,} Princeton Mathematical Series,
 No. 43, Monographs in Harmonic Analysis, III, Princeton University Press, Princeton, NJ, 1993, MR \textbf
{95c}:42002.

\bibitem
{Stein_Weiss} E. M. Stein, G. Weiss, \textit
{Introduction to Fourier Analysis on Euclidean Spaces,} Princeton University Press, Princeton, New Jersey, 1971,
 MR \textbf
{46} \#4102.

\bibitem
{Vega_1988}
 L. Vega, \textit
{Schr\"odinger Equations: Pointwise Convergence to the Initial Data,}
 Proc. Amer. Math. Soc. \textbf
{102} (1988),
 874--878, MR \textbf 
{89d}:35046.

\bibitem
{Tao_Vargas_2000}
 T. Tao, A. Vargas, \textit
{A Bilinear Approach to Cone Multipliers. II. Applications,}
 Geom. Funct. Anal. \textbf
{10} (2000),
 216--258, MR \textbf 
{2002e}:42013.

\bibitem
{Walther_1995}
 B. Walther, \textit
{Maximal Estimates for Oscillatory Integrals with Concave Phase,}
 Harmonic Analysis and Operator Theory (Caracas, 1994), 485--495,
 Contemp. Math., \textbf
{189}, Amer. Math. Soc., Providence, RI, 1995, MR \textbf
{96e}:42024.

\bibitem
{Walther_1999a}
 B. G. Walther, \textit
{Some $L \sp p (L \sp \infty)$- and $L \sp 2 (L \sp 2)$-estimates for Oscillatory Fourier Transforms,}
 Analysis of Divergence (Orono, ME, 1997), 213--231,
 Appl. Numer. Harmon. Anal., Birkh\"aser Boston, Boston, MA, 1999, MR \textbf
{2001e}:42013.

\bibitem
{Walther_1999b}
 B. G. Walther, \textit
{A Sharp Weighted $L \sp 2$-estimate for the Solution to the Time-dependent Schr\"odinger Equation,}
 Ark. mat. \textbf
{37} (1999),
 381--393, MR \textbf
{2000g}:35029.

\bibitem
{Walther_2001}
 B. G. Walther, \textit
{Higher Integrability for Maximal Oscillatory Fourier Integrals,}
 Ann. Acad. Sci. Fenn. Ser. A I Math. \textbf
{26} (2001),
 189--204, MR \textbf
{2002f}:42021.

\bibitem
{Walther_2002a}
 B. G. Walther, \textit
{Estimates with Global Range for Oscillatory Integrals with Concave Phase,}
 Colloq. Math. \textbf
{91} (2002),
157--165, MR \textbf
{2002m}:42006.

\bibitem
{Walther_2002b}
 B. G. Walther, \textit
{Regularity, Decay and Best Constants for Dispersive Equations,}
 J. Funct. Anal. \textbf
{89} (2002),
 325--335, MR \textbf
{2003a}:35218.
 
\bibitem
{Walther_2002c}
 B. G. Walther, \textit
{Sharp Maximal Estimates for Doubly Oscillatory Integral,}
 Proc. Amer. Math. Soc. \textbf
{130} (2002),
 3641--3650, MR \textbf
{2003g}:42033.

\bibitem
{Wang_Si_Lei_1991}
 S. L. Wang, \textit
{On the Weighted Estimate of the Solution Associated with the Schr\"odinger Equation,}
 Proc. Amer. Math. Soc. \textbf
{113} (1991),
 87--92, MR \textbf
{91k}:35066.

\bibitem
{Wang_Sichun_2006}
 S. Wang, \textit
{A Radial Estimate for the Maximal Operator Associated with the Free Schrödinger Equation,}
 Studia Math. \textbf
{176} (2006),
 95-–112, MR \textbf
{2007k}:42056.

\end{thebibliography}
 \end{document}